\magnification=\magstep1
\input amstex

\def\centerbmp#1#2#3{\vskip#2\relax\centerline{\hbox to#1{\special
  {bmp:#3 x=#1, y=#2}\hfil}}}

\def\centereps#1#2#3{\vskip#2\relax\centerline{\hbox to#1{\special
  {eps:#3 x=#1, y=#2}\hfil}}}



\input pictex
\documentstyle{amsppt}
\TagsOnRight
\hsize=5.1in                                                  
\vsize=7.8in
\define\R{{\bold R}}

\vskip 2cm
\topmatter
\dedicatory To the memory of V. A. Rokhlin\enddedicatory

\title Counting zeros of closed 1-forms\endtitle
\rightheadtext{}
\leftheadtext{}
\author  Michael Farber\endauthor
\address
School of Mathematical Sciences,
Tel-Aviv University,
Ramat-Aviv 69978, Israel
\endaddress
\email farber\@math.tau.ac.il
\endemail
\thanks{The research was partially supported by a grant from the 
Israel Academy of Sciences and Humanities and by
the Herman Minkowski Center for Geometry; 
a part of this work was done while the author was
visiting Max-Planck Institut f\"ur Mathematik in Bonn} 
\endthanks
\abstract{This paper suggests new topological lower bounds for 
the number of zeros of closed 1-forms within a given cohomology class.
The main new technical tool 
is the {\it deformation complex}, 
which allows to pass to 
{\it a singular limit} and reduce the original problem with closed 1-form to a
traditional problem with a Morse function.
We show by examples that the present approach may provide stronger estimates 
than the Novikov inequalities.
The technique of the paper also applies to study topology of the set of zeros of
closed 1-forms under Bott non-degeneracy assumptions}
\endabstract
\endtopmatter

\define\C{{\bold C}}
\define\RR{{\Cal R}}

\define\Z{{\bold Z}}

\redefine\P{{\Cal P}}

\define\rk{\operatorname{rk}}

\define\GL{\operatorname{GL}}


\define\rank{\operatorname{rank}}

\def\<{\langle}
\def\>{\rangle}

\define\pd#1#2{\dfrac{\partial#1}{\partial#2}}

\redefine\L{\Cal L}

\def\part{\partial}

\NoBlackBoxes

\def\pp{\frak p}
\def\qq{\frak q}

\def\ml{{\operatorname{Mon_L}}}
\define\calP{\Cal P}

\define\calQ{\Cal Q}

\define\Mat{\operatorname{Mat}}
\def\mte{{\operatorname{Mon_{\tilde E}}}}
\define\disk{\operatorname {disk}}
\define\ind{\operatorname {ind}}

\documentstyle{amsppt}
\nopagenumbers

\heading{\bf \S 1. Line bundles, Dirichlet units, and zeros of closed 1-forms}\endheading
In this section we describe some new topological 
estimates on the number of zeros of closed 1-forms. They use 
homology with coefficients in finite dimensional flat vector bundles (compare \cite{N3}, \cite{P}),
replacing the local systems over the
Novikov ring in the well-known Novikov's theory \cite{N1, N2}. 
We will compute examples (cf. \S 5) 
showing that
the method of this paper 
may produce stronger inequalities 
(and thus predicts existence of a larger number of zeros) than the classical
Novikov's approach. 

Estimates on the numbers of zeros of closed 1-forms found interesting applications in symplectic topology. They were initiated by J.-C. Sikorav \cite{S1, S2}. 
We refer also to the work of H. Hofer and D. Salamon \cite{HS}.

\subheading{1.1}
Let $M$ be a smooth manifold, and let $\omega$ be a closed 1-form on $M$, $d\omega=0$.
A point $p\in M$ is {\it a zero} of $\omega$ if $\omega$ vanishes at this point, i.e. $\omega_p=0$. 
A zero $p\in M$ is called {\it nondegenerate} if $\omega$, 
viewed as a map $M\to T^\ast M$, is transversal to the zero section $M\subset T^\ast M$ of the 
cotangent bundle. As is well-known, this condition is equivalent to the requirement that 
in a neighborhood $U$ of $p$ we may write $\omega=df$, where $f:U\to \R$
is a smooth function and $p$ is a non-degenerated (Morse) critical point of $f$. 
The {\it Morse index }
of $p$ is well defined (as the Morse index of $p$ as the critical point of $f$). 
A closed 1-form $\omega$ is called {\it Morse} if
all its zeros are non-degenerate.

\subheading{1.2} Recall that a complex number $a\in \C^\ast$ is called {\it a Dirichlet
unit} if it is a unit of the ring of integers of an algebraic number field. Equivalently, a Dirichlet unit
is an algebraic integer such that its inverse $a^{-1}$ is also an algebraic integer. 
Dirichlet unit $a\in \C$ is a root of a monic polynomial
$$a^k + \beta_1a^{k-1} + \beta_2a^{k-2} +\dots + \beta_k\, =\, 0\tag1-1$$ 
with $\beta_1, \beta_2, \dots, \beta_{k-1}\in \Z$ and $\beta_k=\pm 1$; this property clearly characterizes the Dirichlet units.

In \cite{F2, F3} I used Dirichlet units to construct an analog of the Lusternik - Schnirel\-man theory
for closed 1-forms (i.e. without assuming that the zeros are all non-degenerate).

Consider a flat complex vector bundle $E$ over $M$. 
It is determined by its monodromy, a linear representation of the fundamental group 
$\pi_1(M,x_0)$ 
on the fiber $E_0$ over the base point $x_0$, given by the parallel transport along loops.
For example, a flat line bundle is determined by a homomorphism $H_1(X;\Z)\to \C^\ast$, where $\C^\ast$ is considered as a multiplicative abelian group.

Let $\xi\in H^1(X;\Z)$ be an integral cohomology class. Given a complex number $a\in \C^\ast$,
we will consider the complex flat line bundle over $M$ with the following property: the monodromy
along any loop $\gamma\in \pi_1(M)$ is the multiplication by $a^{\<\xi,\gamma\>}$. We will denote 
this bundle by $a^\xi$.

{\it A lattice} $\L\subset V$ in a finite dimensional vector space $V$ is a finitely 
generated subgroup with $\rank \L =\dim_\C V$.

We will say that a complex flat bundle $E\to X$ of rank $m$ {\it admits an integral lattice} if its monodromy representation $\pi_1(X,x_0)\to \GL_{\C}(E_0)$ is conjugate to a homomorphism 
$\pi_1(X,x_0)\to \GL_{\Z}(\L_0)$, where $\L_0\subset E_0$ is a lattice in the fiber. This condition
is equivalent to the assumption that $E$ is obtained from a local system $\tilde E$
of finitely generated free 
abelian groups over $X$ by tensoring on $\C$.

\proclaim{1.3. Theorem} Let $M$ be a closed smooth manifold and let $\xi\in H^1(M;\Z)$ 
be an integral cohomology class. Let $E\to M$ be a complex flat bundle admitting an integral lattice.
Let $a\in \C^\ast$ be a complex number, which is not a Dirichlet unit.
Then for any closed 1-form $\omega$ on $M$ having 
Morse type zeros and lying in class $\xi$, the number $c_j(\omega)$ 
of zeros of $\omega$ having index $j$ satisfies
$$c_j(\omega)\, \ge\,  \frac{\dim_\C H_j(M;a^\xi\otimes E)}{\dim E},\qquad j=0,1,2, \dots.
\tag1-2$$ 
Moreover, 
$$\sum_{i=0}^j (-1)^i c_{j-i}(\omega) \, \ge\,  
\sum_{i=0}^j (-1)^i \frac{\dim_\C H_{j-i}(M;a^\xi\otimes E)}{\dim E},\qquad j=0,1,2, \dots.
\tag1-3$$
\endproclaim
A proof of Theorem 1.3 is given in \S 3.
\subheading{Remark 1}
Theorem 1.3 gives interesting estimates already in the simplest case when $E$ is taken to be 
the trivial flat line bundle. In this case, however, inequality (1-3) follows from the Novikov inequalities,
written in the form 
$$\sum_{i=0}^j (-1)^i c_{j-i}(\omega) \, \ge\,  
q_j(\xi) + \sum_{i=0}^j (-1)^i b_{j-i}(\xi),\qquad j=0,1,2, \dots,
\tag1-4$$
cf. \cite{F1}, page 48. 

We will show by examples (cf. \S 5) that by applying Theorem 1.3 with
different flat bundles $E$ one obtains stronger
estimates, than provided by the Novikov inequalities.

\subheading{Remark 2} It is easy to show that for transcendental $a\in \C$ the Betti
number $\dim_\C H_i(M; a^\xi)$ equals the Novikov number $b_i(\xi)$ (and in, particular, 
it is the same for all transcendental $a$. 

\subheading{Remark 3} Consider the function
$$a\in\C^\ast \mapsto \dim_\C H_j(M; a^\xi\otimes E).$$
Then there exist only finitely many numbers 
$a_1, a_2, \dots, a_k \in \C^\ast$ (they are called {\it jump points}) 
so that the corresponding 
Betti number $\dim_\C H_j(M; a^\xi\otimes E)$ is the same for any $a\in \C^\ast$ which is not 
one of the jump points.
Following \cite{BF1}, let us denote by $b_j(\xi; E)$ the value of 
$\dim_\C H_j(M; a^\xi\otimes E)$ for $a$ not a jump point.
The number $b_j(\xi; E)$ is a generalization of the Novikov number $b_j(\xi)$. 
For any of the jump points $a_s$ actually holds
$$\dim_\C H_j(M; a_s^\xi\otimes E) > b_j(\xi; E),$$
i.e. {\it the jumps are always positive.} 

Suppose that the flat bundle $E$ admits an integral lattice.
Then the jump points $a_1, a_2, \dots, a_k $ are 
{\it algebraic numbers} (not necessarily algebraic integers).
If a jump happens at a number, which is not a Dirichlet unit, then Theorem 1.3 applies and 
we obtain estimate (3-2), which is stronger, 
than the inequality 
$$c_j(\omega) \ge \frac{b_j(\xi; E)}{\dim E}.$$

\subheading{Remark 4} {\it Theorem 1.3 becomes false if we allow $a\in \C^\ast$ to be
a Dirichlet unit. }
To explain this, note that any Dirichlet unit $a\in \C$ is an eigenvalue of an integral square 
matrix $B=(b_{ij})$ with $\det(B)=1$. We may find a diffeomorphism of a compact smooth 
manifold $h: F\to F$ so that $h$ induces the matrix $B$ on homology 
of some dimension $k$.
Consider the mapping torus $M$, which obtained from $F\times [0,1]$ by identifying
each point $(x,0)$ with $(h(x),1)$. The manifold $M$ is a smooth fiber bundle over the
circle and so it admits a closed 1-form $\omega$ with no critical points, $c_j(\omega)=0$ for all $j$.
The homology $H_\ast(M; a^\xi)$ is nontrivial if and only if the number $a$ is an eigenvalue
of the monodromy $h_\ast: H_\ast(F;\C)\to H_\ast(F;\C)$. Here $\xi$ denotes the cohomology class
of $\omega$. Hence, if $a$ is a Dirichlet unit, we may construct $M$ so that
$H_\ast(M;a^\xi)\ne 0$ and class
$\xi$ may be realized by a closed 1-form with no critical points.

\heading{\bf \S 2. Morse inequalities for prime ideals}\endheading

In this section we describe Morse type inequalities for prime ideals in commutative rings, which will
be used in the proof of Theorem 1.3. 
The results of this section are known to experts in commutative algebra, 
although I am unable to make a reference.

\subheading{2.1} Let $\RR$ be a commutative Noetherian ring. 
Each prime ideal $\pp\subset \RR$ gives a way of assigning Betti numbers to chain complexes over $\RR$.
Indeed, let $C$ be a chain complex over $\RR$ 
$$0\to C_m \to C_{m-1}\to \dots \to C_0\to 0$$
with finitely generated free $\RR$-modules $C_i$. Given a prime ideal $\pp\in \RR$, we will denote
by $b_i(C,\pp)$ the $\pp$-Betti number of $C$, i.e.
$$b_i(C,\pp)\, =\, \dim_{Q(\RR/\pp)}H_i(C\otimes_\RR Q(\RR/\pp)),\tag2-1$$ 
where $Q(\RR/\pp)$ denotes
the field of fractions of the factor ring $\RR/\pp$. 
Define also {\it the Poincar\'e polynomial} corresponding to $\pp$ as 
$$\P_{C, \pp}(\lambda) = \sum_{i=0}^m \lambda^i b_i(C, \pp).\tag2-2$$

Our purpose is to compare the Poincar\'e polynomials corresponding to two different prime 
ideals $\pp\subset \qq\subset \RR$.

Given two polynomials 
$$\P(\lambda) = p_0+p_1\lambda +\dots + p_m\lambda^m,\quad  
\calQ(\lambda) = q_0+q_1\lambda +\dots + q_{m'}\lambda^{m'}$$ with $p_i, q_i\in \Z$, 
we will write
$$\P(\lambda) \succeq \calQ(\lambda),\tag2-3$$ 
to indicate that the differentce $\P(\lambda)-\calQ(\lambda)$ may be represented in the form 
$$\P(\lambda)-\calQ(\lambda) \, =\, (1+\lambda)\Cal T(\lambda),$$
 where $\Cal T(\lambda)$ is a polynomial with non-negative integral coefficients. 
The relation $\succeq$ defines a partial order on the set of all integral polynomials 
in $\lambda$. 
It is well known that (2-3) is equivalent to 
the following sequence of Morse inequalities:
$$
\sum_{j=0}^r(-1)^j p_{r-j} \ge \sum_{j=0}^r(-1)^j q_{r-j},
\quad  r=0, 1, \dots.
$$

\proclaim{2.2. Theorem} Let $C$ be a free finitely generated chain complex over $\RR$ and let 
$\pp\subset \qq\subset \RR$ be two prime ideals. Then for the corresponding
Poincar\'e polynomials holds
$$\P_{C,\qq}(\lambda) \succeq \calP_{C, \pp}(\lambda).\tag2-4$$ 
\endproclaim
\demo{Proof} Applying twice the Euler - Poincar\'e theorem to the truncated complex
$$0\to C_r \to C_{r-1}\to \dots \to C_0\to 0,$$ 
tensored by $Q(\RR/\pp)$ and by $Q(\RR/\qq)$,
we obtain the following identity
$$
\align
&\sum_{j=0}^{r} (-1)^{j} b_{r-j}(C,\pp) + B_r(\pp) = \\
&\sum_{j=0}^r (-1)^j \rk C_{r-j}=\\
&\sum_{j=0}^r (-1)^j b_{r-j}(C,\qq) + B_r(\qq)
\endalign
$$
where $B_r(\pp)$ denotes the dimension over the field $Q(\RR/\pp)$ of the image of the map
$$d: C_{r+1}\otimes Q(\RR/\pp) \to C_{r}\otimes Q(\RR/\pp);\tag2-5$$
the number $B_r(\qq)$ is defined similarly. 

Our statement now is equivalent to the inequality $B_r(\qq)\le B_r(\pp)$ for all $r$. 

Suppose that the homomorphism 
$d: C_{r+1} \to C_r$ is represented by a matrix with entries in $\RR$. Then $B_r(\qq)$
equals to the maximal size of the minors of this matrix which do not lie in the ideal $\qq$. 
The number $B_r(\pp)$ has the similar descritpion in terms of this matrix. 
Since $\pp\subset \qq$, we obtain $B_r(\qq)\le B_r(\pp)$. \qed
\enddemo

\heading{\bf \S 3. Proof of Theorem 1.3}\endheading

In this section we will describe a proof of Theorem 1.3, 
which will be based on the Morse inequalities
for prime ideals (cf. \S 2) and on the existence of a {\it deformation complex}.  
While proving Theorem 1.3, we will not use the particular construction of the 
deformation complex $C_\ast$;
we will only use its properties (i) - (iv), cf. 3.2 below.

\subheading{3.1} Suppose that we are in conditions of Theorem 1.3.
Since the closed 1-form $\omega$ has an integral cohomology class, 
there exists a map $f: M\to S^1$,
so that $\omega=f^\ast(d\theta)$, where $d\theta$ is the standard angular form on the circle. 
The zeros
of $\omega$ are precisely the critical points of $f$. We assume that $\omega$ has only 
nondegenerated zeros, and hence $f$ is a Morse circle valued function on $M$. 
Choose a regular
value $b\in S^1$ and let $N$ be the result of cutting of $M$ along the codimension one 
submanifold $V=f^{-1}(b)$. We have a canonical identification map $\Pi: N\to M$.
Note that the boundary of $N$ contains two copies of $V$. 
We will denote them by $\partial_+ N$ and $\partial_- N$. 
The notations $\partial_\pm N$ are chosen so that for the normal vector field $X$ on 
$\partial_+N$,
pointing inside $N$,
holds $\Pi^\ast\omega(X)>0$; also, for the normal vector field $Y$ on $\partial_-N$, 
pointing inside $N$, holds $\Pi^\ast\omega(Y)<0$.

\subheading{3.2} In \S 4 we will construct a chain complex $C_\ast$ (which we call {\it deformation
complex}) with the following
properties:
\roster
\item"{(i)}" $C_\ast$ is a free finitely generated chain complex over the polynomial ring 
$P=\Z[\tau]$, where $\tau$ is an indeterminate.
\item"{(ii)}" For any nonzero complex number $a\in \C^\ast$ 
the homology $H_i(\C_a\otimes_P C_\ast)$, as a vector space,
is isomorphic to $H_i(M;a^{-\xi}\otimes E)$, where $i=0, 1, 2, \dots $. 
Here $\C_a$ denotes $\C$ with the $P$-module 
structure, so that $\tau$ acts by multiplication on $a\in \C$.
\item"{(iii)}" Let $\C_0$ denote $\C$ with the trivial $P$-module 
structure, i.e. $\tau$ acts as zero on $\C$. Then
the homology $H_i(\C_0\otimes_P C_\ast)$
is isomorphic to $H_i(N,\partial_+ N;\Pi^\ast E)$, where $i=0, 1, 2, \dots $. 
\item"{(iv)}" Let $p$ be a prime number and let $\Z_p$ denote $\Z/p\Z$, considered as a $P$-module
with the trivial (i.e. $\tau=0$) action of $\tau$. Then the homology 
$H_i(\Z_p\otimes_P C_\ast)$,
as an abelian group, is isomorphic to $H_i(N,\partial_+ N;\Z_p\otimes \Pi^\ast\tilde E)$, 
where $\tilde E$
denotes a local system of free abelian groups on $M$ with $\tilde E\otimes \C = E$.
\endroster 

Intuitively, we may view the deformation complex $C_\ast$ as a polynomial family of 
complexes $C_\ast(M; a^{-\xi}\otimes E)$, where $a\in \C^\ast$ is a parameter, and such that it has
a "singular limit" as $a\to 0$, which is described in (iii) and (iv). 

The deformation complex $C_\ast$ depends on the data $M, \xi, E$ and is not unique.
We will use only existence of $C_\ast$.

\subheading{3.3} Let $a\in \C^\ast$ be a nonzero complex number which is not an 
algebraic integer. Consider the prime
ideal $\pp_a\subset P$
consisting of integral polynomials 
$$q(\tau) = q_0+ q_1\tau+ \dots +q_m\tau^m, \quad q_j\in \Z$$
with $q(a^{-1})=0$. If $C_\ast$ denotes the deformation complex, then, using the property (ii) 
and the notations, introduced in \S 2, we obtain
$$\dim_{Q(P/{\pp_a})} H_i(Q(P/{\pp_a})\otimes_P C_\ast) = 
\dim_C H_i(\C_{a^{-1}}\otimes_P C_\ast)= 
\dim_\C H_i(M; a^{\xi}\otimes E)$$
and therefore
$$\P_{C_\ast, \pp_a}(\lambda) = 
\sum_{i=0}^{\dim M} \lambda^i \dim_\C H_i(M; a^{\xi}\otimes E).\tag3-1$$
Since we assume that $a$ is not an algebraic integer, the ideal $\pp_a$ contains no polynomials
with the free term $q_0$ equals $1$. Hence, we obtain that there exists a prime number $p$,
so that the free terms of all polynomials $q(\lambda)$ lying in $\pp_a$ are divisible by $p$. In other words,
$$\pp_a\, \subset \qq,\quad\text{where}\quad \qq = (p) +(\tau)\subset P;\tag3-2$$ 
in other words, the ideal $\qq$ is generated by $p$ and $\tau$.
Using property (iv) in 3.2, we obtain
$$\dim_{Q(P/{\qq})} H_i(Q(P/{\qq})\otimes_P C_\ast ) = 
\dim_{\Z_p} H_i(N, \partial_+N; \Z_p\otimes \tilde E),$$
and thus
$$\P_{C_\ast, \qq}(\lambda) = 
\sum_{i=0}^{\dim M} \lambda^i \dim_{\Z_p} H_i(N, \partial_+N; \Z_p\otimes \tilde E).\tag3-3$$
By Theorem 2.2 we have 
$$\P_{C_\ast, \qq}(\lambda) \succeq \P_{C_\ast, \pp_a}(\lambda).\tag3-4$$

The form $\Pi^\ast\omega$ is differential of a Morse 
function $g: N\to \R$, $dg= \Pi^\ast\omega$. Function $g$ is constant on each connected
component of $\partial_+N$.  
The critical points of $g$ lie in the
interior of $N$ and they are in 1-1 correspondence with the zeros of $\omega$. 
Hence, the traditional Morse theory, 
applied to the cobordism $N$ and Morse function $g$,  gives 
$$\sum_{i=0}^{\dim M}c_i(\omega)\lambda^i \, \succeq \, (\dim E)^{-1}\cdot
\sum_{i=0}^{\dim M} \lambda^i \dim_{\Z_p} H_i(N, \partial_+N; \Z_p\otimes \tilde E).\tag3-5$$
Now (3-5), combined with (3-3), (3-4), (3-1), gives
$$\sum_{i=0}^{\dim M}c_i(\omega)\lambda^i \, \succeq \, (\dim E)^{-1}\cdot
\sum_{i=0}^{\dim M} \lambda^i \dim_\C H_i(M; a^{\xi}\otimes E),\tag3-6$$
which is equivalent to (1-3). 

We are left to prove Theorem 1.3 assuming that $a^{-1}\in \C^\ast$ is not an algebraic integer.
We will apply the result proven above, where we replace the cohomology class $\xi$ by $-\xi$ 
and the flat vector bundle $E$ by the flat vector
bundle $E^\ast\otimes \frak o_M$, where $\frak o_M$ is the orientation bundle of $M$ (i.e. a flat line
bunlde such that the monodromy along any loop equals $\pm 1$ depending on whether 
the orientation of $M$ is preserved or reversed along the loop). Note that $E^\ast$ as well as 
$\frak o_M$ admit intergral lattices. 

This gives
$$\sum_{i=0}^{\dim M}c_i(-\omega)\lambda^i \, \succeq \, (\dim E)^{-1}\cdot
\sum_{i=0}^{\dim M} \lambda^i \dim_\C H_i(M; a^{-\xi}\otimes E^\ast\otimes \frak o_M).\tag3-7$$
It is clear that $c_i(-\omega) = c_{n-i}(\omega)$, where $n=\dim M$. Also, 
from the Poincar\'e duality we obtain
$\dim_\C H_i(M; a^{-\xi}\otimes E^\ast\otimes \frak o_M)=
\dim_\C H_{n-i}(M; a^{\xi}\otimes E)$. Substitutimg these equalities into (3-7), dividing
both sides by $\lambda^{n}$ and then replacing the indeterminate $\lambda^{-1}$ by
$\lambda$, gives (3-6), which is equivalent to (1-3).

This completes the proof. \qed

Theorem 1.3 may also be deduced from the main result of \cite{FR}.

\heading{\bf \S 4. Construction of the deformation complex}\endheading

In this section we will construct the deformation complex and prove its properties (i) - (iv), which were
used in \S 3. The construction here is slightly different from \cite{F3}, \cite{F4}, although the
obtained complex is essentially the same.

\subheading{4.1. The cell decomposition} In this subsection
we will describe a cell decomposition of $M$, related to the given closed 1-form $\omega$. 
We will repeat some of the constructions of subsection 3.2.

We assume that $\omega$ has an integral indivisible cohomology class. Then
there exists a smooth map
$f: M\to S^1$, so that $\omega=f^\ast(d\theta)$, where $d\theta$ is the angular form on the circle $S^1$.
Clearly $f$ is a Morse map and the zeros of $\omega$ are precisely the critical points of $f$. 
Choose a regular value $b \in S^1$,
and let $N$ be the result of cutting of $M$ along the codimension one 
submanifold $V=f^{-1}(b)$. We have a canonical identification map $\Pi: N\to M$.
The boundary of $N$ contains two copies of $V$. 
We denote them by $\partial_+ N$ and $\partial_- N$. 
The normal vector field $X$ on $\partial_+N$,
pointing inside $N$, satisfies $\Pi^\ast\omega(X)>0$. Also, for the normal vector field $Y$ on $\partial_-N$, pointing into $N$, holds $\Pi^\ast\omega(Y)<0$.

The form $\Pi^\ast\omega$ is exact, i.e. it is differential of a Morse 
function $g: N\to \R$, $dg= \Pi^\ast\omega$. Each connected
component of $\partial_+N$ consists of points of local minimum of $g$.  
The critical points of $g$, which lie in the
interior of $N$, are in 1-1 correspondence with the zeros of $\omega$ and have the 
same indices. Hence $c_i(g)=c_i(\omega)$ for any $i$.

Fix a cell decomposition of $N$ so that $V=\partial_+N$ and $\partial_-N$ be subcomplexes and so
that the natural homeomorphism $J: \partial_+N\to \partial_-N$ be a cellular isomorphism. 

The manifold $M$ is homeomorphic to the factor-space $N/\sim$, 
where we identify each pair of points $v\in \partial_+N$ and $J(v)\in \partial_-N$.
In order to give the factor-space $N/\sim$ a CW-structure, consider a cylinder $V\times [0,1]$
having the standard cell-decomposition: for each $i$-dimensional cell $e$ of $V$ we have 
two $i$-dimensional cells $e\times 0$ and $e\times 1$ and one $(i+1)$-dimensional cell $e\times I$
of the cylinder $V\times I$. Glue each point $(v,1)$ of the top face of the 
cylinder with the point $v\in V\subset N$; also, glue each point $(v,0)$ of bottom face of 
the cylinder with the point $J(v)\in N$. As the result we obtain a 
CW-decomposition 
$$M\, =\, (V\times I\cup N)/\sim. \tag4-1$$

\subheading{4.2. The chain complex}
Here we will calculate the cellular chain complex $C_\ast(\tilde M)$ of the universal covering $\tilde M$
of the CW-complex (4-1). 
It is a complex of free left $\Z\pi$-modules, where $\pi$ is the fundamental group of $M$.

In general the complement $M-V$ may be disconnected. Let us choose base points
$v_1, v_2, \dots, v_l$, one for each connected component of $M-V$. For each $j=1, 2, \dots, l$
we may find a smooth path
$\mu_j:[0,1]\to M$, such that $\mu_j(0)=v_1$, $\mu_j(1)=v_j$ and 
$\mu_j$ has intersection number zero with $V$, i.e. $\mu_j\cdot V=0$. Note that $V$
has a fixed orientation of its normal bundle. To construct $\mu_j$, first connect the points
$v_1$ and $v_j$ by an arbitrary path $\overline \mu_j$ in $M$, and then set 
$\mu_j = \overline \mu_j - (\overline \mu_j\cdot V)\delta$, where $\delta$ is a closed loop in $M$
such that $\delta\cdot V=1$. Such loop $\delta$ exists since we assume that the cohomology class
$\xi=[\omega]\in H^1(M;\Z)$ is indivisible.

Each cell of $X$ can be lifted into the covering $\tilde M$, and all possible lifts are parameterized 
by the elements of the fundamental group $\pi =\pi_1(M, v_1)$.
In order to fix the lifts of the cells of $M$, we will describe for each cell 
$e\subset M$ a path $\gamma_e$ in $M$, starting from the base point $v_1$ and leading to an internal point of $e$. 
We will call $\gamma_e$ {\it the tail of} $e$.
After an arbitrary choice of a lift of the base point $v_1$, we will obtain lifts of all the cells in the covering $\tilde M$, determined (in an obvious way) by the choice of the tails.

For the cells $e\subset N$ we will
choose their tails as follows. Assume that $e$ lies in the component of $N$ containing $v_j$.
Then we set $\gamma_e= \mu_j\sigma_e$, where $\sigma_e$ is a path in $N$ connecting $v_j$
with an internal point of $e$.

The tail of cells of the form $e\times I\subset V$, where $e\subset V$, are constructed as follows.
First travel along the existing tail $\gamma_e$ of cell $e$, which leads from the base point 
$v_1$ to an internal point of 
$e\subset V=V\times 1\subset N$, and then
drop slightly down to an internal point of the cell $e\times I$. 

After the choice of tails as above we obtain a free basis of $C_\ast(\tilde M)$ (over the group ring $\Z\pi$) formed by the lifts of the cells of $M$. 

The boundary homomorphism of $C_\ast(\tilde M)$, 
applied to a cylindrical generator $e\times I$ is given by the formula
$$d(e\times I) = d(e)\times I + (-1)^{i}[e\times 1 - e\times 0], \quad i=\dim e.\tag4-2$$
Here $e\times 1$ can be identified with $e$. The generator $e\times 0\in C_\ast(\tilde M)$
is supported by the cell $J(e)\subset \partial_-N\subset N$, however it has a different tail. Indeed, it is
clear that the tail of the cell $e\times 0$ has intersection number $-1$ with $V$, 
because it first travels along the tail of $e\times I$, and then arrives (staying inside $e\times I$)
at the face $e\times 0 \simeq J(e)$
of $e\times I$. Hence we may rewrite
(4-2) as follows
$$d(e\times I) = d(e)\times I + (-1)^{i}[e - g\cdot J(e)],\tag4-3$$
 where
$$ g\in \pi = \pi_1(M,v_1), \quad \xi(g)=-1.\tag4-4$$
Recall that $\xi\in H^1(M;\Z)$ denotes the cohomology class of $\omega$.
\subheading{4.3. The deformation complex} Let $\Z\pi_-\subset \Z\pi$ denote the subring of the group 
ring generated by the group elements $g\in \pi$ with $\xi(g)\le 0$. 
 From the description of the
chain complex $C_\ast(\tilde M)$ (and in particular, from (4-3), (4-4)) we see that the cells of 
$M$ with their tails specified above, 
generate a free chain complex $C'$ over $\Z\pi_-$ such that
$$\Z\pi\otimes_{\Z\pi_-}C' \, =\, C_\ast(\tilde M).$$
Let $\tilde E\to M$ be a local system of free abelian groups so that $\tilde E\otimes \C = E$; it exists
because we assume that flat bundle $E$ admits an integral lattice. 
We have the monodromy representation
$$\mte: \pi \to \GL(m;\Z),\quad m=\rank(\tilde E).\tag4-5$$
Here for $g,g'\in \pi$ holds $\mte(gg') = \mte(g')\mte(g)$, i.e. (4-5) is an anti-isomorphism.
We obtain a ring anti-homomorphism
$\Z\pi_- \to \Mat(m; P),\quad P=\Z[\tau],$
where for $g\in \pi$ with $\xi(g)\le 0$ we set
$g\mapsto \tau^{-\xi(g)}\mte(g).$
This defines a structure of right $\Z\pi_-$-module on $P^m$, commuting with the obvious 
left $P$-action. Hence, $P^m$ becomes a $(P, \Z\pi_-)$-bimodule.
We may fanally define the deformation complex $C_\ast$ as
$$C_\ast = P^m\otimes_{\Z\pi_-}C'.\tag4-6$$

\subheading{4.4} Now we check that the deformation complex (4-6) satisfies the conditions
(i) - (iv) of subsection 3.2. 
Indeed, (i) is obvious. To show (ii), note that for $a\in \C^\ast$ we have an isomorphism
$$\C_a\otimes_P C_\ast = \C_a\otimes_P (P^m\otimes_{\Z\pi_-}C') \simeq 
(\C_a\otimes_\Z \Z^m)\otimes_{\Z\pi} C_\ast(\tilde M)\tag4-7$$
and $\C_a\otimes_\Z \Z^m\simeq \C^m$ is a right $\Z\pi$-module with respect to the action
$$g\mapsto a^{-\xi(g)} \operatorname{Mon}_E(g)\in \GL(m;\C).$$ 
This clearly implies (ii).

Let us now prove (iv). We have
$$\Z_p\otimes_P C_\ast\simeq \Z_p\otimes_P (P^m\otimes_{\Z\pi_-}C')\simeq 
\Z_p^m\otimes_{\Z\pi_-}C',\tag4-8$$ 
where the right action of $\Z\pi_-$ on $\Z_p^m$ is given by
$$
g\mapsto \cases \mte(g), \quad\text{if}\quad \xi(g)=0\\
0\quad\text{if}\quad \xi(g)<0,
\endcases\tag4-9
$$
where $g\in \pi$. From (4-3), (4-4) we see that complex (4-8) will have as its basis the cells of $N$
and cells of the form $e\times I$, one for each cell $e$ of $V$. The boundary map of (4-8)
acts on the cylindrical generators $e\times I$ as follows
$$d(e\times I) = \partial e\times I + (-1)^{\dim e}\cdot e\tag4-10$$
(the last term in (4-3) disappeares). 
This implies that all the cells $e\subset V$ and $e\times I\subset V\times I$ generate a subcomplex
$\frak N_p$ of the chain complex (4-8). 
Moreover, the map $e\mapsto (-1)^{\dim e}e\times I$ is a chain contraction of $\frak N_p$
(because of (4-10)).

Hence we obtain an isomorphism 
$H_i(\Z_p\otimes_P C_\ast) \simeq H_i((\Z_p\otimes_P C_\ast)/\frak N_p)$ and the homology
of the factor-complex
$(\Z_p\otimes_P C_\ast)/\frak N_p$ is clearly $H_i(N,V;\Z_p\otimes \Pi^\ast\tilde E)$. 
This proves (iv).

The proof of (iii) is similar. We will not use (iii) in this paper. 

\heading{\bf \S 5. An example}\endheading

In this section we show that inequalities (1-2), (1-3) 
produce stronger estimates than the Novikov inequalities. 
The example, which is described here, is a modification of 
example 1.7 in \cite{BF}.

\subheading{5.1}
Let $k_N\subset S^3$ be the connected sum of $N$ copies of the trefoil knot. Let $X$ be closed
3-dimensional manifold obtained by 0-surgery on $k_N$. Then $H^1(X;\Z)\simeq \Z$; we will
denote by $\eta\in H^1(X;\Z)$ a generator. 

Recall that the Alexander polynomial of the trefoil is
$\Delta(\tau) = \tau^2 -\tau +1.$ 
Its two roots we denote by $b, b^{-1}\in \C$; they are
Dirichlet units. 

Consider the rank 2 flat vector bundle $F\to X$, where 
$F\, \, \simeq\, \,  b^\eta \oplus b^{-\eta}.$ 
From the knot theory \cite{R} we know that 
$$\dim_\C H_1(X;F) \, =\, 2N,\quad \dim_\C H_1(X;\C) =0.\tag5-1$$
Let us show that $F$ admits an integral lattice. The monodromy representation of $F$ 
is a 2-dimensional vector space
with a basis $e_1, e_2$, on which the meridian $\gamma\in \pi_1(X)$ of $k_N$ acts as follows
$$e_1\mapsto be_1, \quad e_2\mapsto b^{-1}e_2.$$
Hence the vectors $e_1+e_2$ and $be_1+b^{-1}e_2$ form a basis, 
and the action of $\gamma$ in this basis is represented by an integral matrix
$\left[\matrix 0 & -1\\ 1 & 1\endmatrix\right]$.

\subheading{5.2}
Consider the 3-dimensional manifold given as the connected sum 
$$M=X\# (S^1\times S^2).$$
Thus $M=M_+\cup M_-,$
where $M_+\cap M_-=S^2$ and
$$M_+ = X-\{\disk\}, \quad M_- = (S^1\times S^2) -\{\disk\}.$$

Consider a flat vector bundle $E\to M$ such that $E|_{M_+}\simeq F|_{M_+}$ and 
$E|_{M_-}$ is trivial. Let $\xi\in H^1(M;\Z)$ be a class such that $\xi|_{M_+}=0$ and $\xi|_{M_-}$
is a generator. 

We want to compare $\dim_\C H_1(M;a^\xi)$ with $\dim_\C H_1(M;a^\xi\otimes E)$ 
for different $a\in \C^\ast$. Using the Mayer-Vietoris sequence, we obtain
$$\dim_\C H_1(M;a^\xi) = \dim_\C H_1(M_+;\C)+ \dim_\C H_1(M_-;a^\xi),$$
and hence using (5-1) we obtain that $\dim_\C H_1(M;a^\xi) = 0$ for $a\ne 1$. 
All Novikov numbers $b_i(\xi)=0=q_i(\xi)$ vanish.
However, the same
argument as above, using (5-1) and the Mayer-Vietoris sequence, yields 
$$\dim_\C H_1(M;a^\xi\otimes E) = 2N.$$ 

This example shows that an appropriate choice of flat bundle $E$, appearing in inequalities 
(1-2) and (1-3), may provide stronger estimates, than when taking $E$ to be the trivial bundle.

\heading{\bf \S 6. Closed 1-forms with non-isolated zeros}\endheading

In this section we study the topology of the set of zeros of closed 1-forms under Bott type \cite{B}
nondegeneracy condition. First inequalities of this type were obtained in \cite{BF}.

\subheading{6.1}
Let $\omega$ be a smooth closed real valued 1-form on $M$, $d\omega=0$,
which is assumed to have zeros {\it non-degenerate in the sense of R.Bott}
\cite{B}. This means that the set of points of $C\subset M$, where the form
$\omega$ vanishes, form a submanifold of $M$, and the {\it Hessian} of $\omega$ is
{\it non-degenerate on the normal bundle to} $C$.

If $N$ is a small tubular
neighborhood of $C$ in $M$, then the integral $\int_\gamma\omega$ vanishes along
any loop $\gamma\subset N$ (since $\gamma$ is homologous to a curve in $C$). 
Thus there exists a unique real
valued smooth function $f$ on $N$ such that $df=\omega_{|_N}$ and
$f_{|_C}=0$.
The {\it Hessian} of $\omega$ is then defined as the
Hessian of $f$.

Let $\nu(C)$ denote the normal bundle of $C$ in $M$. Note that
$\nu(C)$ may have different dimension over different connected
components of $C$.  Since Hessian of $\omega$ is non-degenerate, the
bundle $\nu(C)$ splits into the Whitney sum of two subbundles
$$\nu(C)\ =\ \nu^+(C)\ \oplus \nu^-(C),$$
such that the Hessian is strictly positive on $\nu^+(C)$ and strictly
negative on $\nu^-(C)$. Here, the dimension of the bundles
$\nu^+(C)$ and $\nu^-(C)$ over different connected components of $C$
may be different.

For every connected component $Z$ of set $C$, the
dimension of the bundle $\nu^-(C)$ over $Z$ is called the {\it index}
of $Z$ and is denoted by $\ind(Z)$.
Let $o(Z)$ denote
the {\it orientation bundle of $\nu^-(C)_{|_Z}$}, considered as a local system with fiber $\Z$.

The following is a generalization of Theorem 1.3 for closed 1-forms with non-isolated zeros.

\proclaim{6.2. Theorem} Let $\omega$ be a closed 1-form on $M$ having Bott type zeros and 
representing an
integral cohomology class $\xi=[\omega]\in H^1(M;\Z)$. 
Let $\tilde E\to M$ be a local system of free abelian groups.
Let $a\in \C^\ast$ be a complex number, not an algebraic integer; let 
$p$ be a prime number with the property that any integral polynomial $q(\tau)$ with 
$q(a)=0$ has top coefficient divisible by $p$.
Then
$$
\aligned
&\sum_{Z} \sum_{i=0}^n \lambda^{\ind(Z)+i}
\dim_{\Z_p} H_i(Z;\Z_p\otimes \tilde E_{|_Z}\otimes o(Z))
\, \succeq\, \\
&\succeq \, \, \sum_{i=0}^n \lambda^i \dim_\C H_i(M; a^\xi\otimes \tilde E).
\endaligned\tag6-1
$$
In the first sum $Z$ runs over all connected components of the set of zeros 
of $\omega$.
\endproclaim

Note that for transcendental $a$ the prime number $p$ may be taken arbitrarily.
In the case of transcendental $a\in \C^\ast$ inequalities (6-1) turn into the inequalities of \cite{BF}, i.e.
with $\C$ replacing $\Z_p$ in the first sum (6-1).
In fact, the method of the proof, suggested in the present paper, also proves the main result of 
\cite{BF} (by using property (iii) in subsection 3.2 of the deformation complex). 
In \cite{BF} the proof used the method of Witten deformation and analytical tools.

\demo{Proof} We simply repeat all the arguments of the proof of Theorem 1.3 given in \S 3, using the
Morse inequalities for prime ideals and the deformation complex. On the last stage instead of inequality
(3-5) we use now the following inequality
$$
\aligned
&\sum_{Z} \sum_{i=0}^n \lambda^{\ind(Z)+i}
\dim_{\Z_p} H_i(Z;\Z_p\otimes \tilde E_{|_Z}\otimes o(Z))
\, \succeq\, \\
&\succeq \sum_{i=0}^n \lambda^i \dim_{\Z_p} H_i(N, \partial_+N;\Z_p\otimes \tilde E),
\endaligned\tag6-2
$$
which is just a slight generalization of the 
well-known inequality of Bott \cite{B}
with $\Z_p$ coefficients, applied to the manifold with boundary
$N$ and to Bott function $g: N\to \R$, cf. \S 3. \qed
\enddemo

\heading{\bf \S 7. Classes of higher rank}\endheading

In this section we announce a generalization of Theorem 1.3 for cohomology classes $\xi$
of higher rank. 
\subheading{7.1}
Let $M$ be a manifold. We will denote by $H$ the first homology group $H_1(M;\Z)$.
Let $\xi\in H^1(M,\R)$ be a real cohomology class. It can be viewed as a homomorphism
$\xi: H_1(M;\Z)=H\to \R$; we denote by $\ker(\xi)$ the kernel. 
Given a polynomial $p\in \Z[H]$, one defines two numbers
$d_\xi(p)$ ({\it the $\xi$-degree of $p$}) and $v_\xi(p)$ 
({\it the $\xi$-top coefficient}) as follows. Let $p=\sum_{j=1}^n \beta_j h_j$, where 
$\beta_j\in \Z$ and $h_j\in H$. Then  
$d_\xi(p)$ is defined as the maximal number $d=d_\xi(p)\in \R$ such that the sum
$v_\xi(p) = \sum \beta_j$, taken over all $j$ with $\<\xi, h_j\>=d$, is nonzero.

Let $L\to M$ be a complex flat line bundle. We will assume that the monodromy of $L$ is trivial
along any loop in $M$ representing a homology class in $\ker(\xi)$. 
$L$ determines a {\it monodromy homomorphism}
$$\ml : \Z[H]\to \C,\tag7-1$$
assigning to any $h\in H$ the monodromy of $L$
along $h$. We will denote by $\Cal I_L\subset  \Z[H]$ the kernel of the homomorphism $\ml$.

\proclaim{7.2. Definition} (A) 
We will say that a flat complex line bundle $L\to M$ is a {\it $\xi$-algebraic integer}
if (i) the monodromy of $L$ is trivial
along any loop in $M$ representing a homology class in $\ker(\xi)$; and (ii) the ideal $\Cal I_L$
contains a polynomial $p\in \Cal I_L$ with $v_\xi(p)=\pm 1$.

(B) We will say that a complex flat line bundle $L\to M$ is a $\xi$-Dirichlet unit if $L$ and the
dual flat line bundle $L^\ast$ are
$\xi$-algebraic integers.  \endproclaim 

The following statement generalizes Theorem 1.3 on forms with arbitrary cohomology classes.

\proclaim{7.3. Theorem} Let $M$ be a closed smooth manifold and let $\xi\in H^1(M;\R)$ 
be a real cohomology class. Let $E\to M$ be a flat complex vector bundle admitting an integral lattice.
Let $L\to M$ be a flat complex line bundle, which is not a $\xi$-Dirichlet unit.
Then for any closed 1-form $\omega$ on $M$ having 
Morse zeros and lying in the class $\xi$, the number $c_j(\omega)$ of zeros
of $\omega$ having index $j$ satisfies
$$c_j(\omega)\, \ge\,  \frac{\dim_\C H_j(M;L\otimes E)}{\dim E},\qquad j=0,1,2, \dots.\tag7-2$$ 
\endproclaim
The proof will be published elsewhere. 
It uses the Morse inequalities for prime ideals and an analog of the
deformation complex for forms of higher rank $>1$.

\enddocument